\documentclass{amsart}
\usepackage{amsmath,amssymb}
\usepackage{hyperref}
\theoremstyle{plain}
\newtheorem{theorem}{Theorem}[section]
\newtheorem{proposition}[theorem]{Proposition}
\newtheorem{lemma}[theorem]{Lemma}
\newtheorem{maintheorem}{Theorem}

\theoremstyle{remark}
\newtheorem{remark}[theorem]{Remark}
\newtheorem{example}[theorem]{Example}

\newcommand{\Z}{\mathbb Z}
\newcommand{\Q}{\mathbb Q}
\newcommand{\C}{\mathbb C}
\newcommand{\R}{\mathbb R}
\newcommand{\cZ}{\mathcal Z}
\newcommand{\RZ}{\mathbb R \mathcal Z}
\DeclareMathOperator{\lk}{lk}
\newcommand{\wideand}{\quad\text{ and }\quad}

\title[Smooth structures on polyhedral products]{Smooth structures on polyhedral products of disks and spheres}

\author{Suyoung Choi}
\address{Department of Mathematics, Ajou University, 206, World cup-ro, Yeongtong-gu, Suwon 16499, Republic of Korea}
\email{schoi@ajou.ac.kr}

\author{Hyeontae Jang}
\address{School of Mathematics, Korea Institute for Advanced Study,\newline 85 Hoegiro, Dongdaemun-gu, Seoul 02455, Republic of Korea}
\email{hjang1112@kias.re.kr}

\thanks{The first-named author was supported by the National Research Foundation of Korea Grant funded by the Korean Government (RS-2025-00521982). The second-named author was supported by a KIAS Individual Grant (MG105401) at Korea Institute for Advanced Study.}
\date{\today}

\begin{document}

\begin{abstract}
    We prove that every polyhedral product of disks and their boundary spheres that is a topological manifold admits a stably parallelizable smooth structure.
    We also show that if a moment-angle complex or a real moment-angle complex admits a smooth structure for which the coordinate action is smooth, then the link of every face of the underlying simplicial complex is a simplicial sphere, and every such smooth structure is stably parallelizable.
\end{abstract}

\maketitle

\section{Introduction}

The polyhedral product assigns a space to a simplicial complex and a family of pairs of spaces.
Let $K$ be a simplicial complex on the vertex set $[m] = \{1, 2, \dots, m\}$, and let $(\underline{X}, \underline{A}) = \{(X_i, A_i)\}_{i=1}^m$ be a family of pairs of spaces.
The \emph{polyhedral product} of $(\underline{X}, \underline{A})$ over~$K$ is the space
$$
    (\underline{X}, \underline{A})^K = \bigcup_{\sigma \in K} \Bigl( \prod_{i \in \sigma} X_i \times \prod_{i \notin \sigma} A_i \Bigr) \subseteq \prod_{i=1}^m X_i.
$$
We write $(X, A)^K$ when every pair $(X_i, A_i)$ is equal to $(X, A)$.

Polyhedral products were introduced in this generality by Bahri, Bendersky, Cohen and Gitler~\cite{Bahri-Bendersky-Cohen-Gitler2010}.
They unify constructions from toric topology, homotopy theory, geometric group theory, commutative algebra and the theory of subspace arrangements~\cite{Bahri-Bendersky-Cohen2020}.
The \emph{moment-angle complex} $\cZ_K = (D^2, S^1)^K$ originates in the work of Davis and Januszkiewicz on torus actions~\cite{Davis-Januszkiewicz1991}.
It is homotopy equivalent to the complement of an arrangement of coordinate subspaces of $\C^m$~\cite{Buchstaber-Panov2015book}.
It is homeomorphic to an intersection of real quadrics when $K$ is polytopal~\cite{Bosio-Meersseman2006, Gitler-Lopez_de_Medrano2013}.
The moment-angle complex $\cZ_K$ and the \emph{real moment-angle complex} $\RZ_K = (D^1, S^0)^K$ are the central objects of toric topology.
The topology of these spaces and their quotients has been studied in many works (see, for example,~\cite{Buchstaber-Panov2015book, Cai-Choi2021, Choi-Park2017Cohomology, Panov-Veryovkin2016}).

In this paper, we study the polyhedral products that are manifolds.
The basic examples come from disks and their boundary spheres.
Let $d_1, \dots, d_m$ be positive integers, and $n - 1$ the dimension of $K$.
We write $(\underline{D}, \underline{S})^K$ for the polyhedral product of the pairs $(D^{d_i}, S^{d_i - 1})$ over~$K$.
We call $K$ a \emph{generalized homology sphere} if the link of every face $\sigma$ of $K$ has the same integral homology as an $(n-1-|\sigma|)$-sphere.
Here and throughout this paper, the empty set is regarded as a face of $K$ whose link is $K$, and the empty simplicial complex is regarded as a $(-1)$-sphere.
Cai~\cite{Cai2017} proved that $(\underline{D}, \underline{S})^K$ is a topological manifold if and only if $K$ is a generalized homology sphere, where $K$ is also required to be simply connected if every $d_i$ is $1$ and $\dim K \geq 2$.

It is a fundamental problem to decide whether these manifolds admit smooth structures.
When $K$ is polytopal, the intersections of quadrics give smooth structures on $\cZ_K$ and on $\RZ_K$ with stably trivial tangent bundles~\cite{Gitler-Lopez_de_Medrano2013}.
In this case, the stable triviality of the tangent bundles also follows from~\cite[Theorem~6.6]{Davis-Januszkiewicz1991}.
Panov and Ustinovsky constructed a smooth structure on~$\cZ_K$ when $K$ underlies a complete simplicial fan~\cite[Theorem~2.2]{Panov-Ustinovsky2012}.
They asked for which triangulations $K$ of a sphere the manifold $\cZ_K$ admits a smooth structure~\cite[Question~2.3]{Panov-Ustinovsky2012}.
This question is also stated in the book of Buchstaber and Panov \cite[Problem~6.5.3]{Buchstaber-Panov2015book}, which is the standard reference in toric topology.

Our first result gives a complete answer to the above questions.
A smooth manifold is \emph{stably parallelizable} if its tangent bundle is stably trivial.
\begin{maintheorem}\label{thm:main}
    Let $K$ be a simplicial complex on $[m]$, and let $d_1, \dots, d_m$ be positive integers.
    If the polyhedral product $(\underline{D}, \underline{S})^K$ is a topological manifold, then it admits a stably parallelizable smooth structure.
\end{maintheorem}
Since we have the complete criterion for $M=(\underline{D}, \underline{S})^K$ to be a topological manifold, Theorem~\ref{thm:main} indeed determines exactly which $M$ admit smooth structures.
The main point of the proof is that the stable tangent microbundle of $M$ is trivial.
The smooth structure is then given by the smoothing theory of Kirby and Siebenmann~\cite{Kirby-Siebenmann1977book} when the dimension is at least five.
The first step is not special to polyhedral products.
Proposition~\ref{prop:tangential} shows that the stable tangent microbundle of a manifold with an action of a compact Lie group is trivial whenever the fiber inclusion of its Borel fibration is null-homotopic after one suspension.

The torus~$T^m$ acts on~$\cZ_K$ by coordinatewise rotations, and the group~$(\Z_2)^m$ acts on~$\RZ_K$ by coordinatewise reflections.
We call a smooth structure on~$\cZ_K$ or on~$\RZ_K$ \emph{equivariant} if this action is smooth.
The smooth structures constructed in~\cite{Gitler-Lopez_de_Medrano2013} and~\cite{Panov-Ustinovsky2012} are equivariant.
In general, the existence of a smooth structure and the existence of an equivariant smooth structure are different questions.
A \emph{simplicial sphere} is a simplicial complex whose geometric realization is homeomorphic to a sphere.

\begin{maintheorem}\label{thm:equi}
    Let $K$ be a simplicial complex on~$[m]$.
    If $\cZ_K$ or $\RZ_K$ admits an equivariant smooth structure, then the link of every face of $K$ is a simplicial sphere.
    Moreover, every equivariant smooth structure on~$\cZ_K$ or on~$\RZ_K$ is stably parallelizable.
\end{maintheorem}

The converse of the first statement follows from the results of Davis~\cite{Davis2014}.
Hence, $\cZ_K$ admits an equivariant smooth structure if and only if the link of every face of~$K$ is a simplicial sphere, and the same holds for $\RZ_K$.

Throughout this paper, manifolds are closed unless otherwise stated.
Simplicial complexes are finite and have no ghost vertices, that is, every element of $[m]$ is a vertex of a simplicial complex on~$[m]$.

\noindent\textbf{Statement on the Use of AI.}
We used AI to discuss possible approaches to the proofs and to check individual steps. We also used it to identify relevant references, assess their applicability, and improve the English. Both authors have reviewed the arguments and references and take full responsibility for the entire manuscript.

\section{Suspension of the fiber inclusion}

Let $K$ be a simplicial complex on~$[m]$, and let $X_1, \dots, X_m$ be connected CW complexes with basepoints.
We consider the polyhedral product $(\underline{X}, \ast)^K$ of the pairs $(X_i, \ast)$ over~$K$.
Let
$$
    j \colon (\underline{X}, \ast)^K \hookrightarrow \prod_{i=1}^m X_i
$$
be the inclusion.
We write $\Sigma$ for the reduced suspension. 

\begin{lemma}\label{lem:suspension}
    Let $f \colon A \to (\underline{X}, \ast)^K$ be a based map such that $j \circ f$ is based null-homotopic.
    Then $\Sigma f$ is based null-homotopic.
\end{lemma}
\begin{proof}
    The natural suspension splitting of~\cite[Theorems~2.10 and~2.15]{Bahri-Bendersky-Cohen-Gitler2010} gives based homotopy equivalences
    $$
        \Sigma (\underline{X}, \ast)^K \simeq \bigvee_{\varnothing \neq I \in K} \Sigma \Bigl(\bigwedge_{i \in I} X_i \Bigr) \wideand \Sigma \prod_{i=1}^m X_i \simeq \bigvee_{\varnothing \neq I \subseteq [m]} \Sigma \Bigl(\bigwedge_{i \in I} X_i \Bigr).
    $$
    The second one is the case where $K$ is the full simplex on~$[m]$.
    By naturality, $\Sigma j$ corresponds to the inclusion of the wedge summands indexed by the nonempty faces of $K$.
    Collapsing the other summands gives a based map
    $$
        r \colon \Sigma \prod_{i=1}^m X_i \longrightarrow \Sigma (\underline{X}, \ast)^K
    $$
    such that $r \circ \Sigma j$ is based homotopic to the identity.
    Since $j \circ f$ is based null-homotopic, we get
    $$
        \Sigma f \simeq r \circ \Sigma j \circ \Sigma f = r \circ \Sigma(j \circ f) \simeq \ast.
    $$
\end{proof}

Let $G$ be a compact Lie group, and let $M$ be a connected space with a continuous action of $G$.
Fix a universal principal $G$-bundle $EG \to BG$ such that $BG$ is a CW complex.
The \emph{Borel construction} of $M$ is the space $M_G = EG \times_G M$.
Given basepoints $e_0 \in EG$ and $m_0 \in M$, we base $BG$ at $b_0 = [e_0]$ and $M_G$ at $[e_0, m_0]$.
The \emph{Borel fibration} of $M$ is the fiber bundle
$$
    M \xrightarrow{ i } M_G \xrightarrow{ p } BG
$$
where $i(x) = [e_0, x]$ and $p([e, x]) = [e]$.
We call $i$ the \emph{fiber inclusion}.

We consider $\RZ_K$ with the action of $(\Z_2)^m$ by coordinatewise reflections and $\cZ_K$ with the action of $T^m$ by coordinatewise rotations.
We take $EG = (S^\infty)^m$ with $BG = (\R P^\infty)^m$ for $G = (\Z_2)^m$ and $BG = (\C P^\infty)^m$ for $G = T^m$.
Let $j$ denote the inclusion of $(\R P^\infty, \ast)^K$ or $(\C P^\infty, \ast)^K$ into the corresponding $BG$.

By~\cite[Theorem~4.3.2]{Buchstaber-Panov2015book} and its real analogue, there is a homotopy equivalence $h'$ from the domain of $j$ to the corresponding Borel construction such that $j = p \circ h'$ (see also~\cite[Example~2.3]{Davis2012}).
We choose $e_0$ and $m_0$ so that $h'$ maps the basepoint to $[e_0, m_0]$.
We use these choices for $\RZ_K$ and $\cZ_K$ in the rest of the paper.

\begin{lemma}\label{lem:fiber}
    Let $M$ be $\RZ_K$ with $G = (\Z_2)^m$ or let $M$ be $\cZ_K$ with $G = T^m$.
    Then the suspended fiber inclusion $\Sigma i$ is based null-homotopic.
\end{lemma}
\begin{proof}
    Since the basepoints are nondegenerate, $h'$ has a based homotopy inverse~$h$.
    Then $j \circ h$ is based homotopic to $p$.
    Since $p \circ i$ is constant, the map $j \circ h \circ i$ is based null-homotopic.
    By Lemma~\ref{lem:suspension}, the map $\Sigma(h \circ i)$ is based null-homotopic.
    Since $\Sigma h$ is a based homotopy equivalence, the map $\Sigma i$ is based null-homotopic.
\end{proof}

Denote by $[X, Y]_\ast$ the set of based homotopy classes of based maps from~$X$ to~$Y$.

\begin{lemma}\label{lem:loop}
    Let $f \colon A \to X$ be a based map such that $\Sigma f$ is based null-homotopic.
    If $Y$ has the based homotopy type of a loop space, then the map
    $$
        f^\ast \colon [X,Y]_\ast \longrightarrow [A,Y]_\ast
    $$
    that sends $[g]$ to $[g \circ f]$ is constant at the class of the constant map.
\end{lemma}
\begin{proof}
    By assumption, there are a based space $W$ and a based homotopy equivalence $e \colon Y \to \Omega W$.
    Let $g \colon X \to Y$ be a based map, and let $\widetilde g \colon \Sigma X \to W$ be the adjoint of~$e \circ g$ under the adjunction between suspensions and loop spaces.
    Then the adjoint of $e \circ g \circ f$ is
    $$
        \widetilde g \circ \Sigma f \colon \Sigma A \longrightarrow W.
    $$
    This map is based null-homotopic since $\Sigma f$ is.
    Hence, $e \circ g \circ f$ is based null-homotopic, and so is $g \circ f$.
\end{proof}

\section{Smooth structures}

Let $G$ be a compact Lie group, and let $M$ be a connected topological manifold with a continuous action of $G$.
Let $p \colon M_G \to BG$ be the Borel fibration of $M$, and let $i \colon M \to M_G$ be the fiber inclusion.
We use microbundles in the sense of Milnor~\cite{Milnor1964}.
A \emph{microbundle} of dimension $n$ over a space $B$ consists of a space~$E$ and maps $s \colon B \to E$ and $\pi \colon E \to B$ such that $\pi \circ s$ is the identity and the following condition holds: every point of $B$ has an open neighborhood $U$ such that some open neighborhood of $s(U)$ in $\pi^{-1}(U)$ is homeomorphic to $U \times \R^n$ by a homeomorphism that sends $s(x)$ to $(x, 0)$ for all $x \in U$ and commutes with the projections to $U$.
The \emph{tangent microbundle} of $M$ is the diagram
$$
    M \xrightarrow{ \Delta } M \times M \xrightarrow{ \mathrm{pr}_1 } M
$$
where $\Delta$ is the diagonal map and $\mathrm{pr}_1$ is the first projection.
The \emph{vertical tangent microbundle} of $p$ is the diagram
$$
    M_G \xrightarrow{ \Delta } M_G \times_{BG} M_G \xrightarrow{ \mathrm{pr}_1 } M_G.
$$
This diagram is a microbundle since $p$ is a fiber bundle whose fiber is $M$.
Its pullback along $i$ is the tangent microbundle of $M$.

We write $B\mathrm{Top}$ for the classifying space of stable topological bundles.
Let 
$$
    \tau_M^{\mathrm{st}} \in [M, B\mathrm{Top}]_\ast
$$
be the class of a based classifying map of the stable tangent microbundle of $M$.
We write $0$ for the class of the constant map.
The following proposition is a homotopy-theoretic version of the argument in~\cite{Hasui-Kishimoto-Kizu2022}.

\begin{proposition}\label{prop:tangential}
    Let $M$ be a connected topological manifold with a continuous action of a compact Lie group $G$.
    If the suspension~$\Sigma i$ of the fiber inclusion is based null-homotopic, then $\tau_M^{\mathrm{st}} = 0$.
\end{proposition}
\begin{proof}
    The space $M_G$ is paracompact and Hausdorff since $BG$ is a CW complex and $p$ is a fiber bundle with compact fiber.
    By the microbundle representation theorem \cite[Theorem~3.3]{Holm1967}, the vertical tangent microbundle of $p$ is represented by a topological bundle~$\xi$ over~$M_G$.
    Then $i^\ast \xi$ represents the tangent microbundle of~$M$.
    
    Let $g \colon M_G \to B\mathrm{Top}$ be a based classifying map of the stabilization of $\xi$.
    Then $g \circ i$ represents $\tau_M^{\mathrm{st}}$.
    The space $B\mathrm{Top}$ is an infinite loop space \cite[p.~124]{Madsen-Milgram1979book}.
    Since $\Sigma i$ is based null-homotopic, Lemma~\ref{lem:loop} gives
    $$
        \tau_M^{\mathrm{st}} = [g \circ i] = 0.
    $$
\end{proof}

We are now ready to prove our first main theorem.

\begin{theorem}[Theorem~\ref{thm:main}]
    If the polyhedral product $(\underline{D}, \underline{S})^K$ is a topological manifold, then it admits a stably parallelizable smooth structure.
\end{theorem}
\begin{proof}
    By~\cite[Theorem~7.2]{Bahri-Bendersky-Cohen-Gitler2015}, the polyhedral product $(\underline{D}, \underline{S})^K$ is homeomorphic to~$\RZ_L$ for a simplicial complex $L$ (see also~\cite[Lemma~2.7]{Cai2017}).
    The complex $L$ is obtained from $K$ by simplicial wedge operations~\cite{Choi-Park2016Wedge}, and it has no ghost vertices since $K$ has none.
    A stably parallelizable smooth structure on $\RZ_L$ can be pulled back along this homeomorphism.
    Hence, we may assume that every $d_i$ is $1$, that is, $(\underline{D}, \underline{S})^K = \RZ_K$.
    
    Assume that $\RZ_K$ is a topological manifold, and put $M = \RZ_K$.
    Let $n - 1$ be the dimension of $K$.
    Then $M$ is a connected manifold of dimension~$n$, and it is a finite CW complex of dimension~$n$.
    Lemma~\ref{lem:fiber} and Proposition~\ref{prop:tangential} give $\tau_M^{\mathrm{st}} = 0$.
    
    We first consider the case $n \geq 5$.
    Let $B\mathrm{O}$ be the classifying space of stable real vector bundles.
    Since $\tau_M^{\mathrm{st}} = 0$, a classifying map of the stable tangent microbundle of $M$ is null-homotopic.
    Hence, the constant map $M \to B\mathrm{O}$ is a lift of it up to homotopy along the forgetful map $B\mathrm{O} \to B\mathrm{Top}$.
    By the classification theorem for stable smoothings \cite[Essay~IV, Theorem~10.1]{Kirby-Siebenmann1977book}, this lift is realized by a smooth structure on~$M$ whose tangent bundle is stably trivial.
    
    Next, we consider the case $n \leq 4$.
    By~\cite[Theorem~2.3]{Cai2017}, the complex $K$ is a generalized homology sphere, and it is simply connected if $n \geq 3$.
    A generalized homology sphere of dimension at most $2$ is a PL~sphere.
    This applies to $K$ when $n \leq 3$ and to the links of the vertices of $K$ when $n = 4$.
    In the latter case, $K$ is a simply connected PL manifold of dimension~$3$.
    Then $K$ is a PL~sphere by the Poincar\'e conjecture~\cite{Morgan-Tian2007book} and Moise's theorem~\cite{Moise1952}.
    In both cases, $M$ is a PL manifold by~\cite[Theorem~2.3]{Cai2017}, and it admits a compatible smooth structure since its dimension is at most seven~\cite[Part~II, Section~5]{Hirsch-Mazur1974book}.
    
    We show that this smooth structure is stably parallelizable.
    Stiefel--Whitney classes are defined for stable topological bundles.
    Hence, $\tau_M^{\mathrm{st}} = 0$ implies
    $$
    w_1(M) = w_2(M) = 0.
    $$
    Therefore, a classifying map of the stable tangent bundle of $M$ has a lift
    $$
    s \colon M \longrightarrow B\mathrm{Spin}
    $$
    where $B\mathrm{Spin}$ is the classifying space of stable spin vector bundles.
    It is enough to prove that $s$ is null-homotopic.
    The space $B\mathrm{Spin}$ is $3$-connected, and $H^4(B\mathrm{Spin}; \Z)$ is infinite cyclic with a generator $q_1$ such that $2 q_1$ is the pullback of the first Pontryagin class \cite[Section~2.1]{Crowley-Nordstrom2021}.
    Hence, $s$ is null-homotopic if $n \leq 3$.
    Assume that $n = 4$.
    Then $H^4(M; \Z) \cong \Z$ since $M$ is orientable.
    The rational Pontryagin classes depend only on the stable tangent microbundle.
    Hence, the image of $p_1(M)$ in $H^4(M; \Q)$ vanishes, and thus $p_1(M) = 0$.
    Since $2 s^\ast q_1 = p_1(M)$, we get $s^\ast q_1 = 0$.
    Consider a map from $B\mathrm{Spin}$ to the Eilenberg--MacLane space $K(\Z, 4)$ that represents $q_1$.
    Since $\pi_4(B\mathrm{Spin})$ is infinite cyclic~\cite[Section~2.1]{Crowley-Nordstrom2021}, the Hurewicz theorem shows that this map induces an isomorphism on~$\pi_4$, and hence it is $5$-connected.
    Since $M$ is a CW complex of dimension~$4$, two maps $M \to B\mathrm{Spin}$ are homotopic if they pull $q_1$ back to the same class.
    Therefore, $s$ is null-homotopic.
\end{proof}

\begin{remark}
    Let $M$ be a polyhedral product $(\underline{D}, \underline{S})^K$ that is a topological manifold.
    The proof of Theorem~\ref{thm:main} shows that the stable tangent microbundle of $M$ is trivial.
    Hence, every characteristic class of $M$ that depends only on the stable tangent microbundle is trivial.
    In particular, the Stiefel--Whitney classes of $M$ vanish with respect to every smooth structure, which recovers~\cite[Theorems~1.2 and~3.2]{Hasui-Kishimoto-Kizu2022} for $\cZ_K$ and $\RZ_K$, and so do the rational Pontryagin classes.
    Consequently, $M$ is oriented null-cobordant with respect to every smooth structure~\cite{Wall1960}.
\end{remark}

\section{Equivariant smooth structures}

Let $K$ be a simplicial complex on~$[m]$.
We consider the cubical complex
$$
    cc(K) = ([0,1], \{1\})^K
$$
of~\cite[Construction~2.9.11 and Example~4.2.6]{Buchstaber-Panov2015book}.
Taking absolute values of the coordinates identifies both $\cZ_K / T^m$ and $\RZ_K / (\Z_2)^m$ with $cc(K)$.
For a face $\sigma$ of $K$, let $F_\sigma$ be the set of points $x$ of $cc(K)$ such that $x_i = 0$ for all $i \in \sigma$.
Let $F_\sigma^\circ$ be the set of points $x$ of $F_\sigma$ such that $x_i > 0$ for all $i \notin \sigma$.
We call the sets $F_\sigma$ the \emph{faces} of $cc(K)$.
Deleting the coordinates in $\sigma$ identifies $F_\sigma$ with $cc(\lk_K(\sigma))$.
Under this identification, $F_\sigma^\circ$ is the open cone on~$\lk_K(\sigma)$.

\begin{theorem}[Theorem~\ref{thm:equi}]
    If $\cZ_K$ or $\RZ_K$ admits an equivariant smooth structure, then the link of every face of $K$ is a simplicial sphere.
    Moreover, every equivariant smooth structure on~$\cZ_K$ or on~$\RZ_K$ is stably parallelizable.
\end{theorem}
\begin{proof}
    Let $M$ be $\cZ_K$ or $\RZ_K$ with an equivariant smooth structure, and let $G$ be the group that acts on it.
    In the rest of this proof, manifolds need not be compact.
    By~\cite[Theorem~2.3 and Corollary~2.10]{Cai2017}, the complex $K$ is a generalized homology sphere.
    
    Let $\sigma$ be a face of $K$, and let $G_\sigma$ be the coordinate subgroup of $G$ supported on~$\sigma$.
    The fixed point set of $G_\sigma$ in $M$ is a smooth submanifold.
    Let $U_\sigma$ be its open subset where every coordinate outside $\sigma$ is nonzero.
    The compact Lie group~$G / G_\sigma$ acts freely and smoothly on~$U_\sigma$, and hence the orbit space of this action is a smooth manifold.
    The orbit map $M \to cc(K)$ identifies this orbit space with $F_\sigma^\circ$.
    Therefore, $F_\sigma^\circ$ is a topological manifold for every face $\sigma$ of $K$.
    
    Let $L = \lk_K(\sigma)$.
    If $\dim L \leq 1$, then the generalized homology sphere $L$ is a simplicial sphere.
    Assume that $\dim L \geq 2$.
    For a vertex $v$ of $L$, the open star of~$v$ in $|L|$ is the open cone on the link of $v$ in $L$.
    Since this link is $\lk_K(\sigma \cup \{v\})$, the open star is homeomorphic to the manifold $F_{\sigma \cup \{v\}}^\circ$.
    The open stars of the vertices cover~$|L|$, and hence $|L|$ is a topological manifold.
    The open cone $F_\sigma^\circ$ is a contractible manifold of dimension at least $3$.
    By the van Kampen theorem, removing its apex does not change the fundamental group, so the complement is simply connected.
    This complement is homotopy equivalent to~$|L|$.
    Thus, $|L|$ is a simply connected compact manifold with the homology of a sphere.
    By the generalized Poincar\'e conjecture~\cite{Freedman1982, Morgan-Tian2007book, Newman1966}, it is homeomorphic to a sphere.
    This proves the first assertion.
    
    The differential of the action makes the tangent bundle~$TM$ a $G$-vector bundle.
    The Borel construction gives a vector bundle $EG \times_G TM$ over~$M_G$ whose pullback along the fiber inclusion~$i$ is $TM$.
    Let $g \colon M_G \to B\mathrm{O}$ be a based classifying map of the stabilization of this bundle.
    The map $\Sigma i$ is based null-homotopic by Lemma~\ref{lem:fiber}.
    The space $B\mathrm{O}$ has the based homotopy type of a loop space~\cite[p.~124]{Madsen-Milgram1979book}.
    Hence, $g \circ i$ is based null-homotopic by Lemma~\ref{lem:loop}, and thus $TM$ is stably trivial.
\end{proof}

\begin{example}\label{ex:double}
    Let $H$ be a triangulation of the Poincar\'e homology $3$-sphere.
    Let $S^0$ be the simplicial complex that consists of two vertices, and let $K = S^0 \ast S^0 \ast H$ be the join.
    Then $|K|$ is the double suspension of $|H|$, and it is homeomorphic to the $5$-sphere by the double suspension theorem \cite{Cannon1979}.
    Hence, $K$ is a triangulation of a sphere, and both $\cZ_K$ and $\RZ_K$ are topological manifolds.
    By Theorem~\ref{thm:main}, they admit stably parallelizable smooth structures.
    Let $\sigma$ be an edge of $K$ with one vertex in each copy of $S^0$.
    The link of $\sigma$ in $K$ is $H$, and $H$ is not a simplicial sphere.
    By Theorem~\ref{thm:equi}, neither $\cZ_K$ nor $\RZ_K$ admits an equivariant smooth structure.
\end{example}


\providecommand{\bysame}{\leavevmode\hbox to3em{\hrulefill}\thinspace}
\providecommand{\MR}{\relax\ifhmode\unskip\space\fi MR }
\providecommand{\MRhref}[2]{%
  \href{http://www.ams.org/mathscinet-getitem?mr=#1}{#2}
}
\providecommand{\href}[2]{#2}

\end{document}